\documentclass[11pt,a4paper]{article}
\usepackage[margin=1in]{geometry}
\usepackage{amsmath,amsthm,amssymb,enumerate, bbm ,graphicx,color,caption,upgreek, float, tikz,  subcaption,booktabs,longtable, appendix,graphics, pdfpages,rotating,mathtools,textcomp, array, bm}
\DeclarePairedDelimiter\floor{\lfloor}{\rfloor}
\usepackage[pdftex]{hyperref} 
\usepackage[customcolors,shade]{hf-tikz}

\usetikzlibrary{arrows, positioning,decorations.pathreplacing}
\definecolor{blauw}{RGB}{61,158,255}
\definecolor{donkerblauw}{RGB}{0,0,255}
\definecolor{donkergroen}{RGB}{46,148,0}
\definecolor{donkerrood}{RGB}{204,0,0}

\makeatletter 
\newcommand\mynobreakpar{\par\nobreak\@afterheading} 
\makeatother

\newcommand{\N}{\mathbb{N}}
\newcommand{\Z}{\mathbb{Z}}
\newcommand{\C}{\mathbb{C}}
\newcommand{\R}{\mathbb{R}}

\usepackage[english]{babel}

\newtheorem{theorem}{Theorem}[section]

\newtheorem{claim}[theorem]{Claim}
\newtheorem{proposition}[theorem]{Proposition}

\theoremstyle{definition}

\newtheorem*{examp*}{Example}

\makeatletter
\let\@fnsymbol\@arabic
\makeatother
\DeclareMathOperator*{\sgn}{sgn}

\theoremstyle{plain}

\newfloat{Algorithm}{!hbt}{alg}

\newcommand{\ndv}{\textbf{end }}
\newcommand{\ifv}{\textbf{if }}

\newcommand{\printv}{\textbf{print }}

\newcommand{\foreachv}{\textbf{foreach }}

\newcounter{thm}[section]

\title{Semidefinite programming bounds for Lee codes}
\date{}
\author{Sven Polak\thanks{Korteweg-De Vries Institute for Mathematics, University of Amsterdam. E-mail: \href{mailto:s.c.polak@uva.nl}{s.c.polak@uva.nl}. The research leading to these
results has received funding from the European Research Council under the European Union’s Seventh Framework Programme (FP7/2007-2013) / ERC grant agreement \textnumero 339109.}}
\begin{document}
\maketitle
\setcounter{footnote}{1}
\noindent \textbf{Abstract.} For~$q,n,d \in \N$, let~$A_q^L(n,d)$ denote the maximum cardinality of a code~$C \subseteq \Z_q^n$ with minimum Lee distance at least~$d$, where~$\Z_q$ denotes the cyclic group of order~$q$. We consider a semidefinite programming bound based on triples of codewords, which bound can be computed efficiently using symmetry reductions, resulting in several new upper bounds on~$A_q^L(n,d)$.  

The technique also yields an upper bound on the independent set number of the~$n$-th strong product power of the circular graph~$C_{d,q}$, which number is related to the Shannon capacity of~$C_{d,q}$. Here~$C_{d,q}$ is the graph with vertex set~$\Z_q$, in which two vertices are adjacent if and only if their distance (mod~$q$) is strictly less than~$d$. The new bound does not seem to improve significantly over the bound obtained from Lov\'asz theta-function, except for very small~$n$.

\,$\phantom{0}$



\section{Introduction}
  Fix two integers~$n,q \in \N$. Denote by~$\Z_q$ the group of integers~$0,1,\ldots,q-1 \pmod{q}$, which serves as  alphabet. A \emph{word} is an element~$v \in \Z_q^n$ and a \emph{code} is a subset~$C \subseteq \Z_q^n$.  For two words~$u,v \in \Z_q^n$, their \emph{Lee distance} is $d_L(u,v):= \sum_{i=1}^n \min \{ |u_i-v_i|, q-|u_i-v_i|  \}$.   The \emph{minimum Lee distance} $d_{\text{min}}^L(C)$  of a code~$C\subseteq \Z_q^n$ is the minimum of~$d_L(u,v)$ taken over distinct~$u,v \in C$. (If~$|C|\leq 1$, we set~$d_{\text{min}}^L(C)= \infty$.) For any natural number~$d$, define
\begin{align} \label{aleeqnd}
A^L_q(n,d):= \max \{ |C| \, \, | \,\, C \subseteq \Z_q^n, \,\, d_{\text{min}}^L(C) \geq d \}.
\end{align}
 The Lee distance was introduced by C.Y.\ Lee in 1958~\cite{lee}. If~$q=2$ or~$q=3$, the Lee distance coincides with the Hamming distance. 
 For~$q\geq 4$, the Lee distance does not only take into account the number of symbols that are different in two words (which is measured by the Hamming distance), but also to what extent these symbols are different. Because of this property, the Lee distance is used  in certain communication systems for information transmission (so called `phase modulated systems', see  \cite[Chapter 8]{berlekamp}). 
 
Generally, it is an interesting and nontrivial problem to determine~$A_q^L(n,d)$ for given~$q,n,d$. Quistorff made a table of upper bounds on~$A_q^L(n,d)$ based on analytic arguments \cite{quistorff}. H.\ Astola and I.\ Tabus calculated several new upper bounds by linear programming \cite{astola2}, using an adaptation of the classical Delsarte bound based on pairs of codewords~\cite{delsarte} (see also~\cite{astola1}).  
 
For binary codes equipped with the Hamming distance, the Delsarte bound was generalized to a semidefinite programming bound based on triples of codewords by A.\ Schrijver~\cite{schrijver}, and later to a quadruple bound by Gijswijt, Mittelmann and Schrijver~\cite{semidef}. Also, nonbinary codes with the Hamming distance have been considered~\cite{tanaka, onsartikel}, codes with mixed alphabets~\cite{mixed} and constant weight (binary) codes~\cite{cw4,schrijver}. In~$\cite{invariant}$, the authors mention the possibility of applying semidefinite programming to Lee codes and they state that to their best knowledge, such  bounds for Lee codes using triples have not yet been studied.

In this paper, we describe how to efficiently compute a semidefinite programming upper bound~$B_3^L(q,n,d)$ on~$A^L_q(n,d)$ based on triples of codewords, using symmetry reductions, and we calculate this  bound for several values of~$q,n,d$. We only consider~$q \geq 5$, since for~$q=4$, the problem of determining~$A^L_4(n,d)$ is equivalent to determining the maximum size of a binary code of length~$2n$ and minimum distance~$d$ using the Gray map (see, for example,~\cite{gray}). We find several new upper bounds on~$A^L_q(n,d)$, see Table~$\ref{tabellee}$.

%For two instances we construct codes achieving the best known upper bounds and we prove uniqueness of these codes using the semidefinite programming output. 

 \begin{table}[ht] \small
\centering
   \begin{minipage}[t]{.44\linewidth}
   \centering
   \begin{tabular}{| r | r | r|| >{\bfseries}r| r|}
    \hline
    $q$ & $n$ & $d$ &  \multicolumn{1}{>{\raggedright\arraybackslash}b{20mm}|}{\textbf{new upper bound}} & \multicolumn{1}{>{\raggedright\arraybackslash}b{20mm}|}{best upper bound previously known} % &  \multicolumn{1}{>{\raggedright\arraybackslash}b{13mm}|}{Delsarte bound}
    \\\hline 
    5 & 4 & 3  & 62 & $64^l$    \\
    5 & 4 & 4  & 27 & $30^l$    \\
    5 & 4 & 5  & 10 & $11^l$    \\
    5 & 5 & 3  & 270 & $276^l$   \\
    5 & 5 & 5  & 36 & $39^l$    \\
    5 & 5 & 6  & 15 & $18^l$    \\
    5 &6 & 3  & 1170& $1176^l$   \\      
    5 & 6 & 4  & 494 & $520^b$    \\
    5 & 6 & 5  & 149 & $155^l$    \\
    5 & 6 & 6  & 60 & $63^l$    \\                        
    5 & 6 & 7  & 25 & $28^l$    \\                  
         5 & 7 & 3 &  5180 & $5208^{bl}$     \\ 
      5 & 7 & 4  &2183  & $2232^b$     \\                 
    5 & 7 & 5  &590  & $608^l$     \\      
    5 & 7 & 6  & 250  & $284^l$     \\      
    5 & 7 & 7  & 79 & $81^l$     \\              
    5 & 7 & 8  &35  &  $41^l$   \\                   
    \hline
    6 & 3 & 3  & 27 & $29^l$    \\
    6 & 3 & 4  & 14 & $17^l$    \\
     6 & 4 & 4  & 78 &$79^l$     \\        
         6 & 4 & 5  & 22 & $26^l$    \\
    6 & 5 & 3  & 693 &  $699^l$   \\
        6 & 5 & 4  & 366 &  $378^l$   \\ 
             6 & 5 & 5  & 107 &  $114^l$   \\
    \hline 
    \end{tabular}\end{minipage}    \,\,\,\,\,\,
   \begin{minipage}[t]{.44\linewidth}
   \centering
   \begin{tabular}{| r | r | r|| >{\bfseries}r| r|}
    \hline
    $q$ & $n$ & $d$ &  \multicolumn{1}{>{\raggedright\arraybackslash}b{20mm}|}{\textbf{new upper bound}} & \multicolumn{1}{>{\raggedright\arraybackslash}b{20mm}|}{best upper bound previously known} % &  \multicolumn{1}{>{\raggedright\arraybackslash}b{13mm}|}{Delsarte bound}
    \\\hline 
     6 & 5 & 6  & 61 &  $67^l$   \\       
    6 & 5 & 7  & 22 &  $24^{bl}$   \\ 
            6 & 6 & 6 & 273 & $293^l$\\    
        6 & 6 & 7 & 79 & $85^l$\\    
    6 & 6 & 8 & 48 & $52^l$\\    
    6 & 6 & 9 & 16 & $17^l$\\
    \hline 
     7 & 3 & 4  & 21 &$24^{bl}$     \\
    7 & 3 & 5  & 10 &  $11^l$   \\               
     7 & 4 & 3  & 256 & $263^l$   \\
          7 & 4 & 4  & 121 & $128^b$   \\
    7 & 4 & 5 & $\bm{49^*}$ &  $50^l$   \\          
     7 & 4 & 6  & 23 & $27^l$    \\
     7 & 4 & 7  & 11 & $13^l$    \\
      7 & 4 & 8  & 6 &  $7^{bl}$   \\
           7 & 5 &  3 & 1499 &  $1512^l$   \\         
     7 & 5 &  4 &  686&  $720^b$   \\         
     7 & 5 &  5 & 240 &  $249^l$   \\         
     7 & 5 &  6 & 116  &  $130^l$   \\         
     7 & 5 &  7 & 49 &  $54^l$   \\      
     7 & 5 &  8 & 25 &  $28^l$   \\           
          7 & 5 &  9 & 13 & $14^l$   \\ 
                       7 & 6 & 10 & 26 & $31^l$  \\
                  7 & 6 & 11 & 13 & $14^b$   \\   
                  &&&&\\
                  \hline 
    \end{tabular}\end{minipage}    
  \caption{\label{tabellee}\small An overview of the new upper bounds for Lee codes. The new upper bounds are instances of the bound~$B_3^L(q,n,d)$ from~$(\ref{Bleend})$ below. The superscript $^l$ refers to a bound obtained by Astola and Tabus  using linear programming~\cite{astola2}. The superscript~$^b$ refers to a bound from Quistorff~\cite{quistorff}. The superscript~$^*$ refers to an upper bound matching the known lower bound:~$A_7^L(4,5)=49$ is achieved by a linear code~\cite{astolalinear}. }
\end{table}  
 
In Section~\ref{circbounds}, we show how to adapt the new bound to an upper bound~$B_3^{L_{\infty}}(q,n,d)$ on the independent set number of the~$n$-th strong product power of the circular graph~$C_{d,q}$, which number is related to the Shannon capacity of~$C_{d,q}$. The \emph{circular graph}~$C_{d,q}$ is the graph with vertex set~$\Z_q$, in which two vertices are adjacent if and only if their distance (mod~$q$) is strictly less than~$d$. The new bound does not seem to improve significantly over the bound obtained from Lov\'asz theta-function, except for very small~$n$.

\subsection{The semidefinite programming bound}
We define a hierarchy of semidefinite programming upper bounds on~$A^L_q(n,d)$, which is an adaptation of the semidefinite programming hierarchy for binary codes defined by Gijswijt, Mittelmann and Schrijver in~$\cite{semidef}$. For~$k \in \Z_{\geq 0}$, let~$\mathcal{C}_k$ be the collection of codes~$C \subseteq \Z_q^n$ with~$|C|\leq k$.  For any~$D \in \mathcal{C}_k$, we define
\begin{align} 
\mathcal{C}_k(D) := \{C \in \mathcal{C}_k \,\, | \,\, C \supseteq D, \, |D|+2|C\setminus D| \leq k  \}.  
\end{align} 
Note that with this definition~$|C \cup C'| \leq k$ for all~$C,C' \subseteq \mathcal{C}_k(D)$.   
Also we define, for any function~$x : \mathcal{C}_k \to \R$ and~$D \in \mathcal{C}_k$, the~$\mathcal{C}_k(D) \times \mathcal{C}_k(D)$-matrix~$M_{k,D}(x)$ by $M_{k,D}(x)_{C,C'} : = x(C \cup C')$, for~$C,C' \in \mathcal{C}_k(D)$. Now define the following number:
\begin{align} \label{Bleend}
B^L_k(q,n,d):=    \max \{ \sum_{v \in \Z_q^n} x(\{v\})\,\, &|\,\,x:\mathcal{C}_k \to \R, \,\, x(\emptyset )=1, \,\, x(S)=0 \text{ if~$d_{\text{min}}^L(S)<d$}, \notag \\[-1.1em]
& \,\, M_{k,D}(x) \text{ is positive semidefinite for each~$D$ in~$\mathcal{C}_k$}\}. 
\end{align}
 \begin{proposition} \label{trivleeprop}
Fix~$k\in \N$. For all~$q,n,d \in \N$, we have~$A^L_q(n,d) \leq B^L_k(q,n,d)$.
\end{proposition}
\proof
Let~$C \subseteq \Z_q^n$ be a code with~$d_{\text{min}}^L(C)\geq d$ and~$|C| = A_q^L(n,d)$. Define~$x: \mathcal{C}_k \to \R$ by~$x(S)=1$ if~$S \subseteq C$ and~$x(S)=0$ else. Then~$x$ satisfies the conditions in~$(\ref{Bleend})$, where the last condition is satisfied since~$M_{k,D}(x)_{C,C'}=x(C)x(C')$ for all~$C,C'\in \mathcal{C}_k(D)$. Moreover, the objective value equals~$\sum_{v \in  \Z_q^n} x(\{v\}) =|C|=A^L_q(n,d) $, which gives~$ B^L_k(q,n,d) \geq A^L_q(n,d) $.
\endproof
 It can be shown that the bound~$B_2^L(q,n,d)$ is equal to the Delsarte bound in the Lee scheme, which was calculated for many instances by Astola and Tabus in~$\cite{astola2}$. In this paper we consider the bound~$B_3^L(q,n,d)$. The method for obtaining a symmetry reduction, using representation theory of the dihedral and symmetric groups, is an adaptation of the method in~$\cite{onsartikel}$.

 \subsection{Symmetry reductions}
 Fix~$k \in \N$.  Let~$D_q$ be the dihedral group of order~$2q$ and let~$S_n$ be the symmetric group on~$n$ elements. The group~$H:=D_q^n \rtimes  S_n $ acts naturally on~$\mathcal{C}_k$, and this action maintains minimum distances and cardinalities of codes~$C \in \mathcal{C}_k$. We can assume that the optimum~$x$ in~$(\ref{Bleend})$  is $H$-invariant, i.e., $g \circ x = x$ for all~$g \in H$. Indeed, if~$x$ is any optimum solution for~$(\ref{Bleend})$, then for each~$g \in H$, the function~$g \circ x $ is again an optimum solution, since the objective value of~$g \circ x $ equals the objective value of~$x$ and~$g \circ x $ still satisfies all constraints in~$(\ref{Bleend})$. Since the feasible region is convex, the optimum~$x$ can be replaced by the average of~$g \circ x $ over all~$g \in H$. This gives an $H$-invariant optimum solution.

 Let~$\Omega_k$ be the set of~$H$-orbits on~$\mathcal{C}_k$. Then~$|\Omega_k|$ is bounded by a polynomial in~$n$, for fixed~$q$. Since there exists an~$H$-invariant optimum solution, we can replace, for each~$\omega \in \Omega_k$ and~$C \in \omega$, each variable~$x(C)$ by a variable~$z(\omega)$. Hence, the matrices~$M_{k,D}(x)$ become matrices~$M_{k,D}(z)$ and we have considerably reduced the number of variables in~$(\ref{Bleend})$.

We only have to check positive semidefiniteness of~$M_{k,D}(z)$ for one code~$D$ in each~$H$-orbit of~$\mathcal{C}_k$, as for each~$g \in H$, the matrix~$M_{k,g(D)}(z)$ can be obtained by simultaneously permuting rows and columns of~$M_{k,D}(z)$. 

We sketch how to reduce these matrices in size. For~$D \in \mathcal{C}_k$, let~$H_D$ be the subgroup of~$H$ consisting of all~$g \in H$ with~$g(D)=D$. Then the action of~$H$ on~$\mathcal{C}_k$ induces an action of $H_D$ on~$\mathcal{C}_k(D)$.  The simultaneous action of~$H_D$ on the rows and columns of~$M_{k,D}(z)$ leaves~$M_{k,D}(z)$ invariant. This means that the matrices~$M_{k,D}(z)$ are elements of  $(\C^{\mathcal{C}_k(D) \times \mathcal{C}_k(D)})^{H_D}$, which is naturally isomorphic to the \emph{centralizer algebra} of the action of~$H_D$ on~$\C^{\mathcal{C}_k(D)}$, i.e., the collection of~$H_D$-equivariant endomorphisms~$\C^{\mathcal{C}_k(D)} \to \C^{\mathcal{C}_k(D)} $. Therefore, there exists a block-diagonalization~$M_{k,D}(z) \mapsto U^T  M_{k,D}(z) U$ of~$M_{k,D}(z)$, for a matrix~$U$ depending on~$H_D$ but not depending on~$z$. Then~$M_{k,D}(z)$ is positive semidefinite if and only if each of the blocks is positive semidefinite. There are several equal (or equivalent) blocks and after removing duplicate (or equivalent) blocks we obtain a matrix of order bounded polynomially in~$n$, for fixed~$q$, where the entries in each block are linear functions in the variables~$z(\omega)$ (with coefficients bounded polynomially in~$n$). Hence, we have reduced the size of the matrices involved in our semidefinite program.

The reductions of the optimization problem will be described in detail in Section~$\ref{reduct}$. Table~$\ref{tabellee}$ contains the new upper bounds. All improvements have been found using multiple precision versions of SDPA~\cite{nakata}.
 
\section{Preliminaries on representation theory}
In this section we give the definitions and  notation from representation theory (mostly concerning  the symmetric group) used throughout the paper, similarly to the notation used in~$\cite{onsartikel}$. Proofs are omitted, but for more information, the reader can consult Sagan~$\cite{sagan}$. The content of this section is  the same as Section 2 of~\cite{onsartikel, cw4}, so readers who are familiar with one of these papers can safely skip this section. 

A \emph{group action} of a group~$G$ on a set~$X$ is a group homomorphism~$\phi: G \to S_X$, where~$S_X$ is the group of bijections of~$X$ to itself. If~$G$ acts on~$X$, we write~$g \circ x:= \phi(g)(x)$ for all~$g \in G$ and~$x \in X$ and we write~$X^G$ for the set of elements of~$X$ invariant under the action of~$G$. If~$X$ is a linear space, the elements of~$S_X$ are assumed to be linear functions. The action of~$G$ on a set~$X$   induces an action of~$G$ on the linear space~$\mathbb{C}^X$, by~$(g \circ f)(x):= f(g^{-1} \circ x)$, for~$g\in G$,~$f \in \mathbb{C}^X$ and~$x \in X$.    

If~$m \in \N$ and~$G$ is a finite group acting on~$V=\C^m$, then~$V$ is a~$G$\emph{-module}. If~$V$ and~$W$ are~$G$-modules, then a $G$\emph{-homomorphism} (or: $G$\emph{-equivariant map}) $\psi: V \to W$ is a linear map such that~$g \circ \psi(v)=\psi(g \circ v)$ for all~$g \in G$,~$v \in V$. Moreover, a module~$V$ is called \emph{irreducible} if the only~$G$-invariant submodules of~$V$ are~$\{0\}$ and~$V$ itself.

Suppose that~$G$ is a finite group acting \emph{unitarily} on~$V=\C^m$. This means that for each~$g\in G$ there is a unitary matrix~$U_g \in \C^{m \times m}$ such that~$g \circ x = U_gx$ for all~$x \in \C^m$. Consider the inner product~$\langle x, y \rangle :=x^*y$ for~$x,y \in \C^m$, where~$x^*$ denotes the conjugate transpose of~$x \in \C^m$. Then~$V$ can be decomposed as a direct sum of~$G$\emph{-isotypical components}~$V_1,\ldots,V_k$. This means that~$V_i$ and~$V_j$ are orthogonal for distinct~$i$ and~$j$ (with respect to the mentioned inner product), and each~$V_i$ is a direct sum~$V_{i,1} \oplus \ldots \oplus V_{i,m_i}$ of irreducible and mutually isomorphic~$G$-modules, such that~$V_{i,j}$ and~$V_{i',j'}$ are isomorphic if and only if~$i=i'$.  

For each~$i \leq k$ and~$j \leq m_i$ we choose a nonzero vector~$u_{i,j} \in V_{i,j}$ with the property that for each~$i$ and all~$j,j'\leq m_i$ there exists a~$G$-isomorphism~$V_{i,j} \to V_{i,j'}$ mapping~$u_{i,j}$ to~$u_{i,j'}$. For each~$i \leq k$, we define~$U_i$ to be the matrix~$[u_{i,1},\ldots ,u_{i,m_i}]$ with columns~$u_{i,j}$ ($j=1,\ldots,m_i$). Any set of matrices~$\{U_1,\ldots,U_k\}$ obtained in this way is called a \emph{representative} set for the action of~$G$ on~$\C^m$. Then the map
\begin{align} \label{PhiC}
    \Phi : (\C^{m \times m})^G \to \bigoplus_{i=1}^k \C^{m_i \times m_i} \,\, \text{ with } \,\, A \mapsto \bigoplus_{i=1}^k U_i^* A U_i
\end{align}
is bijective. So~$ \dim((\C^{m \times m})^G) =\sum_{i=1}^k m_i^2$, which can be considerably smaller than~$m$. Another crucial property for our purposes is that any~$A \in (\C^{m \times m})^G $ is positive semidefinite (i.e., self-adjoint with all eigenvalues nonnegative) if and only if the image~$\Phi(A)$ is positive semidefinite, i.e., each of the matrices~$U_i^* A U_i$ is positive semidefinite.

It turns out that all representative sets we define consist of real matrices. Then
\begin{align} \label{PhiR} 
\Phi(A) = \bigoplus_{i=1}^k U_i^T A U_i \text{ for }A \in (\R^{m \times m})^G, \,\,\,\,\text{ and }\Phi\left((\R^{m \times m})^G\right) = \bigoplus_{i=1}^k \R^{m_i \times m_i}.
\end{align}
Also,~$A \in \R^{ m \times m}$ is positive semidefinite if and only if each of the matrices~$U_i^T A U_i$ is positive semidefinite ($i=1,\ldots,k$). This is very useful for checking whether~$A$ is positive semidefinite. 

It is convenient to note that, since~$V_{i,j}$ is the linear space spanned by~$G \circ u_{i,j}$ (for each~$i,j$), we have
\begin{align} \label{Rm}
\C^m = \bigoplus_{i=1}^k\bigoplus_{j=1}^{m_i} \C G \circ u_{i,j},
\end{align}
where~$\C G$ denotes the group algebra of~$G$. It will also be convenient to consider the columns of~$U_i$ as elements of the dual space~$(\C^m)^*$ via the inner product mentioned above.

\subsection{A representative set for the action of~\texorpdfstring{$S_n$}{Sn} on~\texorpdfstring{$V^{\otimes n}$}{Vtensorn}}\label{repr}
Fix~$n\in \N$ and a finite-dimensional vector space~$V$. We will consider the natural action of~$S_n$ on~$V^{\otimes n}$ by permuting the indices. We describe a representative set for the action of~$S_n$ on~$V^{\otimes n}$ that will be used repeatedly in the reductions throughout this paper. 

A \emph{partition}~$\lambda $ of~$n$ is a sequence~$(\lambda_1,\ldots, \lambda_h)$ of natural numbers with~$\lambda_1 \geq \ldots \geq \lambda_h >0$ and~$\lambda_1 + \ldots + \lambda_h = n$. The number~$h$ is called the \emph{height} of~$\lambda$. We write~$\lambda \vdash n$ if~$\lambda$ is a partition of~$n$. The~\emph{Young shape} (or \emph{Ferrers diagram})~$Y(\lambda)$ of~$\lambda$ is the set
\begin{align}
    Y(\lambda) := \{(i,j) \in \N^2 \, | \, 1 \leq j \leq h, \, 1 \leq i \leq \lambda_j\}.  
\end{align}
Fixing an index~$j_0 \leq h$, the set of elements~$(i,j_0)$ (for~$1 \leq i \leq \lambda_j$) in~$Y(\lambda)$ is called the~$j_0$\emph{-th row} of~$Y(\lambda)$. Similarly, fixing an element~$i_0 \leq \lambda_1$, the set of elements~$(i_0,j)$ (where~$j$ varies) in~$Y(\lambda)$ is called the~$i_0$\emph{-th column} of~$Y(\lambda)$.  Then the \emph{row stabilizer}~$R_{\lambda}$ of~$\lambda$ is the group of permutations~$\pi$ of~$Y(\lambda)$ with~$\pi(Z)=Z$ for each row~$Z$ of~$Y(\lambda)$. Similarly, the \emph{column stablizer}~$C_{\lambda}$ of~$\lambda$ is the group of permutations~$\pi$ of~$Y(\lambda)$ with~$\pi(Z)=Z$ for each column~$Z$ of~$Y(\lambda)$. 

A \emph{Young tableau} with shape~$\lambda$ (also called a~$\lambda$\emph{-tableau}) is a function~$\tau : Y(\lambda) \to \N$. A Young tableau with shape~$\lambda$ is \emph{semistandard} if the entries are nondecreasing in each row and strictly increasing in each column. Let~$T_{\lambda,m}$ be the collection of semistandard $\lambda$-tableaux with entries in~$[m]$. Then~$T_{\lambda,m} \neq \emptyset$  if and only if~$m$ is at least the height of~$\lambda$. We write~$\tau \sim \tau'$ for~$\lambda$-tableaux~$\tau,\tau´$ if~$\tau'=\tau r$ for some~$r \in R_{\lambda}$.

Let~$B=(B(1),\ldots,B(m))$ be an ordered basis of~$V^*$. For any~$\tau \in T_{\lambda,m}$, define 
\begin{align} \label{utau}
    u_{\tau,B}:= \sum_{\tau'\sim \tau } \sum_{c \in C_{\lambda}} \text{sgn}(c) \bigotimes_{y \in Y(\lambda)} B\left(\tau '(c(y))\right).
\end{align}
Here the Young shape~$Y(\lambda)$ is ordered by concatenating its rows. Then (cf.~$\cite{sagan}$ and~$\cite{onsartikel}$) the set 
\begin{align} \label{reprsetdef}
\left\{\, [u_{\tau,B} \,\, | \,\, \tau \in  T_{\lambda,m}] \,\, | \,\, \lambda \vdash n  \right\},
\end{align}
consisting of matrices, is a representative set for the natural action of~$S_n$ on~$V^{\otimes n}$, for any ordering of the elements in~$T_{\lambda,m}$.

\subsection{A representative set for the action of~\texorpdfstring{$G^n \rtimes S_n$}{G power n rtimes Sn} on~\texorpdfstring{$V^{\otimes n}$}{Vtensorn}}\label{reprmult}
Let~$G$ be any group acting unitarily on~$V:=\C^m$. Suppose that a representative set for the action of~$G$ on~$\C^m$ is given. Here each~$B_i$ is an~$m \times m_i$ matrix, for given integers~$k,m_1,\ldots,m_k$.

Let ${\bm{N}}$ be the collection of all $k$-tuples $(n_1,\ldots,n_k)$ of nonnegative integers adding up
to $n$.
For $\bm{n}=(n_1,\ldots,n_k)\in{\bm N}$, let
$\bm{\lambda\vdash n}$ mean that $\bm{\lambda}=(\lambda_1,\ldots,\lambda_k)$ with
$\lambda_i\vdash n_i$ for $i=1,\ldots,k$.
(So each $\lambda_i$ is equal to a partition $(\lambda_{i,1},\ldots,\lambda_{i,t})$ of~$n_i$, for some $t$.)

For $\bm{\lambda}\vdash\bm{n}$ define
$$
W_{\bm{\lambda}}:=T_{\lambda_1,m_1}\times\cdots\times T_{\lambda_k,m_k},
$$
and for $\bm{\tau}=(\tau_1,\ldots,\tau_k)\in W_{\bm{\lambda}}$ define
\begin{align}\label{16de15d}
v_{\bm{\tau}}:=\bigotimes_{i=1}^ku_{\tau_i,B_i}.
\end{align}
 Proposition~$2$ of~\cite{onsartikel} implies the following.  (In reference~$\cite{onsartikel}$, it is stated that~$G=S_q$ and~$V=\R^{q \times q}$, but with a straightforward adaptation one obtains the following result.)
\begin{proposition}\label{prop2}
The matrix set
\begin{align}\label{matset}
\{~~
[v_{\bm{\tau}}
\mid
\bm{\tau}\in W_{\bm{\lambda}}]
~~\mid
\bm{n}\in\bm{N},\bm{\lambda\vdash n}
\}
\end{align}
is representative for the action of $H:=G^n\rtimes S_n$ on $V^{\otimes n}$ (for any ordering of the elements in~$W_{\lambda}$).
\end{proposition}
Note that the representative set from~$\eqref{matset}$ is real if we start with a real representative set~$\{B_1,\ldots,B_k\}$ for the action of~$G$ on~$V$.

\section{Reduction of the optimization problem \label{reduct}}

In this section we give the reduction of optimization problem~$(\ref{Bleend})$ for computing~$B^L_3(q,n,d)$, using the representation theory from the previous section. First we consider block diagonalizing~$M_{3,D}(z)$ for~$D \in \mathcal{C}_3$ with~$|D|=1$.  Subsequently we consider the case~$D=\emptyset$. Note that for the cases~$|D|=2$ and~$|D|=3$ the matrix~$M_{3,D}(z) = (z(D))$ has order~$1 \times 1$, so it is its own block diagonalization. Hence, in those cases, $M_{3,D}(z)$ is positive semidefinite if and only if~$z(D) \geq 0$. 

\subsection{The case~\texorpdfstring{$|D|=1$}{|D|=1}\label{D1}}
The Lee isometry group~$H=D_q^n \rtimes S_n$ acts transitively on~$\Z_q^n$, so we may assume that~$D=\{\bm{0}$\}, where~$\bm{0}=0\ldots0$ is the all-zero word. The rows and columns of~$M_{3,D}(z)$ are indexed by sets of the form~$\{\bm{0}, \alpha\}$ for~$\alpha \in \Z_q^n$. Then the subgroup~$H_D$ of~$H$ that leaves~$D$ invariant is equal to  $S_2^n \rtimes  S_n$, as the zero word must remain fixed (so we cannot apply a rotation of the alphabet in any coordinate position). Here the non-identity element of~$S_2$ acts on~$\Z_q$, where we consider~$0,\ldots,q-1$ as vertices of a regular~$q$-gon, as a \emph{reflection} switching vertices~$i$ and~$q-i$ (for~$i=1,\ldots,\floor{\frac{q-1}{2}}$). So vertex~$0$ is fixed if~$q$ is odd, and vertices~$0$ and~$q/2$ are fixed if~$q$ is even. For~$i=0,\ldots,q-1$, let~$e_i$ be the~$i$th unit vector of~$\C^{\Z_q}$.
\begin{proposition}
 A representative matrix set for the reflection action of~$S_2$ on~$\C^{\Z_q}$ is
\begin{align} \label{oddrepr}
\{B_1,B_2\}, \,\,\,\text{ with } \,\,
B_1:=\left[e_0, \left( e_{i} + e_{q-i} \right)_{i=1}^{\floor{\frac{q}{2}}}\right], \,\,\,\, B_2:= \left[ \left(e_{i}- e_{q-i}\right)_{i=1}^{\floor{\frac{q-1}{2}}}\right].
\end{align}
\end{proposition}
\proof
  For~$j=1,\ldots,\floor{q/2}+1$, define~$W_{1,j}$ to be the~$1$-dimensional vector space spanned by the~$j$th column~$w_{1,j}$ of~$B_1$. Moreover, for~$j=1,\ldots,\floor{(q-1)/2}$, define~$W_{2,j}$ to be the~$1$-dimensional vector space spanned by the~$j$th column~$w_{2,j}$ of~$B_2$.  Note that each~$W_{i,j}$ is~$S_2$-stable and that~$W_{i,j}$ and~$W_{i',j'}$ are orthogonal whenever~$(i,j) \neq (i',j')$ (with respect to the inner product $u,v \mapsto v^*u$). Observe that, for~$j,j'$ and~$l,l'$ the maps~$W_{1,j}\to W_{1,j'}$ and~$W_{2,l} \to W_{2,l' }$ defined by~$w_{1,j}\mapsto w_{1,j'}$ and~$w_{2,l}\mapsto w_{2,l'}$, respectively, are~$S_2$-equivariant.  Note that the number of~$W_{1,j}$ we have defined is $\floor{q/2}+1$, the number of~$W_{2,j}$ is $\floor{(q-1)/2}$, and
\begin{align*}
(\floor{q/2}+1)^2 + \floor{(q-1)/2}^2 = \begin{cases}
\mbox{$\frac{1}{2}$}q^2+2 &  \text{if~$q$ is even}, \\
(q^2+1)/2  & \text{if~$q$ is odd}, 
\end{cases}
\end{align*}
which is equal to~$|(\Z_q \times \Z_q)/S_2| = \dim ( \C^{\Z_q} \otimes \C^{\Z_q})^{S_2}$. (If~$q$ is even, the points~$(0,0)$, $(q/2,0)$, $(0,q/2)$ and $(q/2,q/2)$ in $\Z_q \times \Z_q$ are fixed by the nonidentity element in~$S_2$. If~$q$ is odd, only the point~$(0,0)$ in~$\Z_q \times \Z_q$ is fixed by the nonidentity element in~$S_2$.)  It follows that the~$W_{1,j}$ and~$W_{2,j}$ form a \emph{decomposition} of~$\C^{\Z_q}$ into \emph{irreducible representations} (as any further representation, or decomposition, or equivalence among the~$W_{i,j}$ would yield that the sum of the squares of the multiplicities of the irreducible representations is strictly larger than~$\dim ( \C^{\Z_q} \otimes \C^{\Z_q})^{S_2}$, which contradicts the fact that~$\Phi$ in~$\eqref{PhiC}$ is bijective). So the matrix set~$\eqref{oddrepr}$ is indeed representative for the action of~$S_2$ on~$\C^{\Z_q}$.  
\endproof 
Note that the representative set is real. Set~$m_1:=\floor{q/2}+1$ and~$m_2:=\floor{(q-1)/2}$. 
Let ${\bm{N}}$ be the collection of all $2$-tuples $(n_1,n_2)$ of nonnegative integers adding up
to $n$. As before, for $\bm{n}=(n_1,n_2)\in{\bm N}$, let
$\bm{\lambda\vdash n}$ mean that $\bm{\lambda}=(\lambda_1,\lambda_2)$ with
$\lambda_i\vdash n_i$ for $i=1,2$.
(So each $\lambda_i$ is equal to a partition $(\lambda_{i,1},\ldots,\lambda_{i,t})$ of~$n_i$, for some $t$.)

For $\bm{\lambda}\vdash\bm{n}$ define
\begin{align} \label{wlambda}
W_{\bm{\lambda}}:=T_{\lambda_1,m_1}\times T_{\lambda_2,m_2},
\end{align}
and for $\bm{\tau}=(\tau_1,\tau_2)\in W_{\bm{\lambda}}$ define
\begin{align}\label{tau1}
v_{\bm{\tau}}:=u_{\tau_1,B_1} \otimes u_{\tau_2,B_2}.
\end{align}
Then Proposition~\ref{prop2} implies that
\begin{align}\label{reprsetdq1}
\{~~
[v_{\bm{\tau}}
\mid
\bm{\tau}\in W_{\bm{\lambda}}]
~~\mid
\bm{n}\in\bm{N},\bm{\lambda\vdash n}
\}
\end{align}
is representative for the action of $S_2^n\rtimes S_n$ on $(\C^{\Z_q})^{\otimes n} =\C^{\Z_q^n}$. Note that the representative set is real.

\subsubsection{Computations for~\texorpdfstring{$|D|=1$ }{|D|=1}} \label{d1comp}
Let~$D = \{ \bm{0} \} \in \mathcal{C}_3$ and let~$\Omega_3$ denote the set of all~$D_q^n \rtimes S_n$-orbits of codes in~$\mathcal{C}_3$. For each~$\omega \in \Omega_3$, we define the~$\mathcal{C}_3(D) \times \mathcal{C}_3(D)$-matrix~$N_{\omega}$ with entries in~$\{0,1\}$ by
\begin{align}
    (N_{\omega})_{\{\bm{0},\alpha\},\{\bm{0},\beta\}} := \begin{cases} 1 &\mbox{if } \{\bm{0},\alpha,\beta\} \in \omega,  \\ 
0 & \mbox{else.} \end{cases} 
\end{align}
Given~$\bm{n}= (n_1,n_2) \in \bm{N}$, for each~$\bm{\lambda} \vdash \bm{n}$ we write~$U_{\bm{\lambda}}$ for the matrix in~$(\ref{reprsetdq1})$ that corresponds with~$\bm{\lambda}$.  For each~$z : \Omega_3 \to \mathbb{R}$ we obtain with~$(\ref{PhiR})$ that
\begin{align} \label{blocks1tt}
    \Phi(M_{3,D}(z)) = \Phi \left(\sum_{\omega \in \Omega_3}z(\omega) N_{\omega} \right ) = \bigoplus_{\bm{n} \in \bm{N}} \bigoplus_{\bm{\lambda} \vdash \bm{n}} \sum_{\omega \in \Omega_3}z(\omega) U_{\bm{\lambda}}^T N_{\omega} U_{\bm{\lambda}}. 
\end{align}
The number of~$\bm{n} \in \bm{N}$, $\bm{\lambda} \vdash \bm{n}$, and the numbers~$|W_{\bm{\lambda}}|$ and~$|\Omega_3|$  are all bounded by a polynomial in~$n$. This implies that the number of blocks in~$(\ref{blocks1tt})$, the size of each block and the number of variables occurring in all blocks are polynomially bounded in~$n$. We now show how to compute the entries of the matrix~$U_{\bm{\lambda}}^T N_{\omega} U_{\bm{\lambda}}$, for all~$\omega \in \Omega_3$,~$\bm{n} \in \bm{N}$, $\bm{\lambda} \vdash \bm{n}$, in polynomial time. That is, we show how to compute the coefficients~$v_{\bm{\tau}}^T N_{\omega}  v_{\bm{\sigma}}$, for~$\bm{\tau}, \bm{\sigma} \in W_{\bm{\lambda}}$, in the blocks~$\sum_{\omega \in \Omega_3} z(\omega)U_{\bm{\lambda}}^T N_{\omega} U_{\bm{\lambda}}$ in polynomial time. 

Let~$\Pi$ be the set of those words that appear as lexicographically minimal element in a~$D_q$-orbit of~$\Z_q^3$. So there is a bijection between~$\Pi$ and the set of orbits of the action of~$D_q$ on~$\Z_q^3$. For any word~$v \in \Z_q^3$, write~$\pi(v)$ for the element in~$\Pi$ that is in the same~$D_q$-orbit of~$\Z_q^3$ as~$v$. Note that
\begin{align}
    \Pi = \{ 00j  \,\, | \,\, j=0,\ldots, \floor{q/2}\} \cup \{ 0jh  \,\, | \,\, j=1,\ldots, \floor{q/2}, h=0,\ldots,q-1\}.
\end{align}
For any element~$P \in \Pi$, define
\begin{align}
    d_P := \sum_{\substack{i,j \in \Z_q: \\ \pi(0ij)=P}} e_{i} \otimes e_{j}.
\end{align}
Then the set~$Z:=\{ d_P \,\, | \,\, P \in \Pi \}$ forms a basis for~$(\C^{\Z_q} \otimes \C^{\Z_q})^{S_2}$, where we consider the reflection action of~$S_2$ on~$\Z_q$, i.e., we consider~$0,\ldots,q-1 \in \Z_q$ as vertices of a regular~$q$-gon, and the non-identity element of~$S_2$ switches the vertices~$i$ and~$q-i$ (for~$i=1,\ldots,\floor{\frac{q-1}{2}}$). We write~$Z^*$ for the dual basis.

Let~$Q$ denote the set of monomials of degree~$n$ on~$(\C^{\Z_q} \otimes \C^{\Z_q})^{S_2}$. Then the function~$(\Z_q^n)^3 \to \mathcal{C}_3$ that maps an ordered triple~$(\alpha,\beta,\gamma)$ to the unordered triple $\{\alpha,\beta,\gamma\}$ induces a surjective function
$
r\, \, : \,\, Q \to \Omega_3 \setminus \{ \{ \emptyset \}\}.
$
For any~$\mu \in Q$, define
$$
K_{\mu}:=\sum_{\substack{d_1,\ldots,d_n\in Z \\ d_1^*\cdots d_n^*=\mu}}\bigotimes_{j=1}^nd_j.
$$
Then a routine calculation (as in Lemma~1 of \cite{onsartikel}) implies that, for each $\omega\in\Omega_3$,
$$
 N_{\omega}=\sum_{\substack{\mu\in Q\\ r(\mu)=\omega}}K_{\mu}.
$$
For any~$\bm{\tau}, \bm{\sigma} \in W_{\bm{\lambda}}$, define the following degree~$n$ polynomial on~$ (\C^{\Z_q} \otimes \C^{\Z_q})^{S_2}$:
\begin{align}\label{polylee}
p_{\bm{\tau,\sigma}}:=
\prod_{i=1}^2
\sum_{\substack{\tau_i'\sim\tau_i \\ \sigma_i'\sim\sigma_i}}\sum_{c_i,c_i'\in C_{\lambda_i}}\sgn(c_ic_i')
\prod_{y\in Y(\lambda_i)}
B_i(\tau_i'c_i(y))\otimes B_i(\sigma_i'c_i'(y)).
\end{align}
This polynomial can be computed (i.e., expressed as linear combination of monomials
in $B_i(j)\otimes B_i(h)$) in time polynomially bounded in $n$, for fixed~$q$ (cf.~\cite{gijswijt,onsartikel}).
%, a description of Gijswijt's method in the above language can also be found in~$\cite{cw4}$
Then a straightforward calculation, highly similar to the one in Lemma~2 of~\cite{onsartikel}, yields that
\begin{align} \label{kmubelangrijk}
\sum_{\mu \in Q} (v_{\bm{\tau}}^T K_{\mu}v_{\bm{\sigma}}) \mu =p_{\bm{\tau,\sigma}}.
\end{align} 
So $\sum_{\mu \in Q} v_{\bm{\tau}}^T K_{\mu}v_{\bm{\sigma}} \mu$  can be computed by
expressing the polynomial $p_{\bm{\tau,\sigma}}$ as linear combination of monomials $\mu\in Q$,
which are products of linear functions in $Z^*$. In order to express $p_{\bm{\tau,\sigma}}$ as linear combination of monomials $\mu\in Q$ it remains to express each $B_i(j)\otimes B_i(h)$
as a linear function into the basis $Z^*$, that is, to
calculate the numbers $(B_i(j)\otimes B_i(h))(d_P)$ for all $i=1,2$ and
$j,h=1,\ldots,m_i$, and $P\in\Pi$. We find
\begin{align} 
B_1(1) \otimes B_1(1) &= 1 d_{000}^*, \notag   \\ 
B_1(1) \otimes B_1(j+1) &= 2 d_{00j}^*,    \text{ for~$j=1,\ldots,\floor{q/2}$} \notag \\
B_1(j+1) \otimes B_1(1) &= 2 d_{0j0}^*,    \text{ for~$j=1,\ldots,\floor{q/2}$}  \notag\\
B_1(j+1) \otimes B_1(h+1) &= 2 d_{0jh}^* + 2 d_{0j(q-h)}^*,    \text{ for~$j,h \in \{1,\ldots,\floor{q/2}\}$}, \notag  \\ 
B_2(j) \otimes B_2(h) &= 2 d_{0jh}^* - 2 d_{0j(q-h)}^*,    \text{ for~$j,h \in \{1,\ldots,\floor{(q-1)/2}\}$},   \label{formuleslee1}
\end{align}  
where the
coefficient of~$d_P^*$ is obtained by evaluating $(B_i(j)\otimes B_i(h))(d_P)$. 
Now one computes the entry $\sum_{\omega \in \Omega_3} z(\omega) \bm{v_{\bm{\tau}}}^T   N_{\omega} \bm{v_{\bm{\sigma}}}$ by first expressing~$p_{\bm{\tau}, \bm{\sigma}}$ as a linear combination of~$\mu \in Q$ and subsequently replacing each $\mu \in Q$ in~$p_{\bm{\tau}, \bm{\sigma}}$ with the variable~$z(r(\mu))$.

\subsection{The case~\texorpdfstring{$D=\emptyset$}{D=empty}\label{Dempty}}
 Let~$D=\emptyset$. The rows and columns of~$M_{3,\emptyset}(z)=M_{2,\emptyset}(z)$ are indexed by words in~$\Z_q^n$ together with the empty set, and~$H_D$ is equal to  $D_q^n \rtimes  S_n$. Here~$D_q$ acts on~$\C^{\Z_q}$  by permuting the vertices~$0,\ldots,q-1$ of a regular~$q$-gon. To compute the block diagonalization of~$M_{2,\emptyset}(z)$, one can use the Delsarte formulas in the Lee scheme~\cite{astola1,astola2}. Here we give the reduction in terms of representative sets.
 
  Let~$\zeta = e^{2 \pi i /q}$ be a primitive~$q$th root of unity.
For each~$j=0,\ldots,\floor{q/2}$, define the vectors $a_j := (1,\zeta^j, \zeta^{2j}, \ldots, \zeta^{(q-1)j})^T$, $b_j :=  (1, \zeta^{-j}, \zeta^{-2j},\ldots, \zeta^{-(q-1)j})^T  \in \C^{\Z_q}$ and set~$V_j := \text{span}\{a_j,b_j\}$. Furthermore, put
$$
c_j:=\frac{\sqrt{\dim V_j} }{2} (a_j+b_j)= \sqrt{\dim V_j}(1, \cos(2j\pi/q), \ldots,\cos(2(q-1)j\pi/q))^T \in \R^{\Z_q} \subseteq \C^{\Z_q}.
$$
 \begin{proposition}  
A representative set for the action of~$D_q$ on~$\C^{\Z_q}$ is given by
\begin{align} \label{linrepr}
\left\{C_1,\ldots,C_{\floor{\frac{q}{2}}+1}\right\}, \,\,\, \text{where } C_j:= c_{j-1}, \text{ for } j=1,\ldots, \floor{\frac{q}{2}}+1.
\end{align}
\end{proposition}
\proof 
Observe that each~$V_j$ is~$D_q$-stable and that~$c_j \in V_j$. Moreover,~$V_l$ and~$V_j$ are orthogonal if~$l \neq j$ (with respect to the inner product $u,v \mapsto u^*v$). To see this, note that~$x:=\zeta^{\pm j \pm l}$ is a~$q$th root of unity unequal to~$1$ if~$j \neq l \in \{ 0,\ldots ,\floor{ q/2 }\}$, so~$1+x+x^2+\ldots+x^{q-1}=0$. This implies that $a_j^*a_l=b_j^*a_l=a_j^*b_l=b_j^*b_l=0$, so~$V_l$ and~$V_j$ are orthogonal. Note that~$\sum_{j=0}^{\floor{q/2}} 1^2 =  \floor{q/2}+1$, which is the number of distinct~$V_j$, is equal to the dimension of~$(\C^{\Z_q \times \Z_q})^{D_q}$. So the~$V_j$ form an orthogonal \emph{decomposition} of~$\C^{\Z_q}$ into \emph{irreducible} representations (as any further representation, or decomposition, or equivalence among the~$V_j$ would yield that the sum of the squares of the multiplicities of the irreducible representations is strictly larger than~$\floor{q/2}+1$, which contradicts the fact that~$\Phi$ in~$\eqref{PhiC}$ is bijective). As~$C_{j+1}$ is an element of~$V_j$ for~$j=0,\ldots,\floor{q/2}$, this implies that~$\left\{C_1,\ldots,C_{\floor{\frac{q}{2}}+1}\right\} $ is a representative matrix set.
\endproof 
Note that the representative set is real, and that each~$C_i$ is a~$q \times 1$-matrix. For simplicity of notation, set~$s:=\floor{\frac{q}{2}}+1$. Let ${\bm{M}}$ be the collection of all $s$-tuples $(n_1,\ldots,n_s)$ of nonnegative integers adding up
to $n$. For~${\bm{n}} \in {\bm{M}}$, define~$v_{\bm{n}}:= C_1^{\otimes n_1} \otimes C_2^{\otimes n_2} \otimes \ldots \otimes C_{s}^{\otimes n_s}$. Proposition~\ref{reprmult} gives the following.
\begin{proposition} 
The set
\begin{align}\label{reprsetdq2}
\{~
v_{\bm{n}}
~\mid
\bm{n}\in{\bm{M}}
\}
\end{align}
is representative for the action of $D_q^n\rtimes S_n$ on $(\C^{\Z_q})^{\otimes n} = \C^{\Z_q^n} = \C^{\mathcal{C}_3(\emptyset) \setminus \{ \emptyset \}}$.
\end{proposition}
Observe that~$D_q^n \rtimes S_n$ acts trivially on~$\emptyset$. The~$D_q^n \rtimes S_n$-isotypical component of~$\mathbb{C}^{\Z_q^n}$ consisting of the~$D_q^n \rtimes S_n$-invariant elements corresponds to the matrix in the representative set indexed by~$\bm{n} = (n,0,\ldots,0)$. Hence  we add a new unit base vector~$\epsilon_{\emptyset}$ to this matrix (as a column) in order to obtain a representative set for the action of~$D_q^n \rtimes S_n$ on~$ \C^{\Z_q^n \cup \{\emptyset\}} = \C^{\mathcal{C}_3(\emptyset)}$. 

\subsubsection{Computations for~\texorpdfstring{$D=\emptyset$ }{D=empty}} \label{d0comp}

In this section we explain how to compute the coefficients in the block diagonalization of~$M_{2,\emptyset}(z)$. First we give a reduction of~$M_{2,\emptyset}(z)$ without the row and column indexed the empty code. Later we explain how the empty code is added. For each~$\omega \in \Omega_2$, we define the~$\Z_q^n\times \Z_q^n$-matrix~$N_{\omega}'$ with entries in~$\{0,1\}$ by
\begin{align}
    (N_{\omega}')_{\alpha,\beta} := \begin{cases} 1 &\mbox{if } \{\alpha,\beta\} \in \omega,  \\ 
0 & \mbox{else.} \end{cases} 
\end{align}
 For each~$z : \Omega_2 \to \mathbb{R}$ we obtain with equations~$(\ref{reprsetdq2})$ and~$(\ref{PhiR})$ that $\Phi \left(\sum_{\omega \in \Omega_2}z(\omega) N_{\omega}' \right ) = \bigoplus_{\bm{n} \in \bm{M}} \sum_{\omega \in \Omega_2} z(\omega) v_{\bm{n}}^T N_{\omega}' v_{\bm{n}}$. 
This shows that~$\Phi \left(\sum_{\omega \in \Omega_2}z(\omega) N_{\omega}' \right)$ is a diagonal matrix. Note that~$|\Omega_2|$ and~$|\bm{M}|$ are polynomially bounded in~$n$.  Now we show how to compute~$v_{\bm{n}}^T N_{\omega}'  v_{\bm{n}}$, for~$\bm{n} \in \bm{M}$ in polynomial time.

Define~$    \Pi' = \{ 00j  \,\, | \,\, j=0,\ldots, \floor{q/2}\} \subseteq \Pi$. For any element~$P \in \Pi'$, define
\begin{align}
    f_P := \sum_{\substack{i,j \in \Z_q: \\ \pi(iij)=P}} e_{i} \otimes e_{j}
\end{align}
Then the set~$\tilde{Z}:=\{ f_P \,\, | \,\, P \in \Pi' \}$ forms a basis for~$(\C^{\Z_q} \otimes \C^{\Z_q})^{D_q}$. Let~$\tilde{Z}^*$ denote the dual basis.
Let~$Q'$ denote the set of monomials of degree~$n$ on~$(\C^{\Z_q} \otimes \C^{\Z_q})^{D_q}$. The function~$(\Z_q^n)^2 \to \mathcal{C}_2$ that maps~$(\alpha,\beta)$ to $\{\alpha,\beta\}$ induces a surjective function $r' \, : \, Q' \to \Omega_2 \setminus \{ \{ \emptyset \}\}$.
For any~$\bm{n} \in \bm{M}$, define the following degree~$n$ polynomial on~$ (\C^{\Z_q} \otimes \C^{\Z_q})^{D_q}$:
\begin{align}\label{polylee2}
p_{\bm{n}}:=
\prod_{i=1}^s (C_i\otimes C_i)^{n_i}.
\end{align}
For any~$\mu \in Q'$, define
$$
K_{\mu}':=\sum_{\substack{f_1,\ldots,f_n\in Z \\ f_1^*\cdots f_n^*=\mu}}\bigotimes_{j=1}^nf_j, \,\,\, \text{ so that (cf. Lemma~2 of~\cite{onsartikel}) }\sum_{\mu \in Q'} v_{\bm{n}}^T K_{\mu}'v_{\bm{n}} \mu =p_{\bm{n}}.
$$
So $\sum_{\mu \in Q'} v_{\bm{n}}^T K_{\mu}'v_{\bm{n}} \mu$ can be computed by
expressing the polynomial $p_{\bm{n}}$ as linear combination of monomials $\mu\in Q'$,
which are products of linear functions in $\widetilde{Z}^*$. To this end, we express each $C_i\otimes C_i$
as linear function into the basis $\widetilde{Z}^*$, i.e., we
calculate the numbers $(C_i\otimes C_i)(f_P)$ for all $i=1, \ldots,s$
and $P\in\Pi'$. We find
\begin{align} 
\text{for $q$ even: } \,\, C_i \otimes C_i &= q \left(f_{000}^*+(-1)^i f_{00(q/2)}^* + 2\sum_{j=1}^{q/2-1} \cos(2 \pi ji/q)f_{00j}^* \right) ,  \text{ for~$i \in \{0,\ldots,q/2\}$}, \notag \\
\text{for $q$ odd: } \,\, C_i \otimes C_i &= q \left(f_{000}^* + 2\sum_{j=1}^{(q-1)/2} \cos(2 \pi ji/q)f_{00j}^*\right) ,  \text{ for~$i \in \{0,\ldots,(q-1)/2\}$}.
 \label{formuleslee2compact}
\end{align} 
Now, as $N_{\omega}'=\sum_{\substack{\mu\in Q'\\\ r'(\mu)=\omega}}K_{\mu}'$ (for each $\omega\in\Omega_2$), one computes the entry $\sum_{\omega \in \Omega_2} z(\omega) v_{\bm{n}}^T   N_{\omega}' v_{\bm{n}}$ by first expressing~$p_{\bm{n}}$ as a linear combination of~$\mu \in Q'$ and subsequently replacing each $\mu \in Q'$ in~$p_{\bm{n}}$  with the variable~$z(r'(\mu))$.

To add the empty code, we add an extra row and column corresponding to the vector~$\epsilon_{\emptyset}$ to the matrix in the representative set indexed by~$\bm{n} = (n,0,\ldots,0)$, as explained below Proposition~\ref{reprsetdq2}. So the only matrix block affected by the empty code in the block diagonalization of~$M_{2,\emptyset}(z)$ is
\begin{align} \label{emptycodeblock}
T:= [\epsilon_{\emptyset}, v_{\bm{n}}]^T M_{2,\emptyset}(z) [\epsilon_{\emptyset}, v_{\bm{n}}], 
\end{align}
 which is a~$2 \times 2$-matrix.   Then~$\epsilon_{\emptyset}^T M_{2,\emptyset}(z)\epsilon_{\emptyset} =M_{2,\emptyset}(z)_{\emptyset,\emptyset}= x(\emptyset)=1$ by definition, see~$(\ref{Bleend})$. Since~$v_{\bm{n}}=C_1^{\otimes n}$ is the all-ones vector, we have $ \epsilon_{\emptyset}^T M_{\emptyset}(z)   v_{\bm{n}} = q^n z_{\omega_0}$,
where~$\omega_0 \in \Omega_{2}$ is the (unique) $D_q^n \rtimes S_n$-orbit of a code of size~$1$.

\section{The strong product power of circular graphs\label{circbounds}}

For any graph~$G=(V,E)$, let~$G^n$ denote the graph with vertex set~$V^n$ and edges between two distinct vertices~$(u_1,\ldots,u_n)$ and~$(v_1,\ldots,v_n)$ if and only if for all~$i \in \{1,\ldots,n\}$ one has either $u_i=v_i$ or~$u_iv_i \in E$. The \emph{Shannon capacity} of~$G$ is defined as
\begin{align}
    \Theta(G):= \sup_{d \in \N} \sqrt[d]{\alpha(G^d)},
\end{align}
where~$\alpha(G^d)$ denotes the maximum cardinality of an independent set in~$G^d$, i.e., a set of vertices no two of which are adjacent~$\cite{shannon}$.

For two integers~$d,q$ with~$q \geq 2d$, the \emph{circular graph} $C_{d,q}$ is the graph with vertex set~$\Z_q$, the cyclic group of order~$q$, in which two distinct vertices are adjacent if and only if their distance (mod~$q$) is strictly less than~$d$. So~$C_{2,q} = C_q$, the circuit on~$q$ vertices.  A classical upper bound on~$\alpha(C_{d,q}^n)$ is given by Lov\'asz's  $\vartheta$-function (see~\cite{lovasz}): one has~
\begin{align}
\alpha(C_{d,q}^n)     \leq \vartheta(C_{d,q}^n)=\vartheta(C_{d,q})^n.
\end{align}
Hence, Lov\'asz's $\vartheta$-function gives an upper bound on the Shannon capacity of~$C_{d,q}$. A closed formula for~$\vartheta(C_{d,q})$ is given in~$\cite{bachoc}$.  We describe how the bound~$B_3^L(q,n,d)$ can be adapted to an upper bound~$B_3^{L_{\infty}}(q,n,d)$ on~$\alpha(C_{d,q}^n)$, which either improves or is equal to the bound obtained from  Lov\'asz's $\vartheta$-function. However, the new bound is not multiplicative over the strong product, so it does not give an upper bound on~$\Theta(C_{d,q})$. 

For distinct~$u,v$ in~$\Z_q^n$, define their \emph{Lee\textsubscript{$\infty$}-distance} $d_{L_{\infty}}(u,v)$ to be the maximum over the distances of~$u_i$ and~$v_i$ (mod~$q$), where~$i$ ranges from~$1$ to~$n$. The \emph{minimum Lee\textsubscript{$\infty$}-distance}~$d_{\text{min}}^{L_{\infty}}(D)$ of a set~$D \subseteq \Z_q^n$ is the minimum Lee\textsubscript{$\infty$}-distance  between any pair of distinct elements of~$D$. (If~$|D|\leq 1$, set~$d_{\text{min}}^{L_{\infty}}(D)=\infty$.) Then~$d_{\text{min}}^{L_{\infty}}(D) \geq d$ if and only if~$D$ is independent in~$C_{d,q}^n$. Define, for~$k \geq 2$,
\begin{align} \label{Bleemaxnd}
B^{L_{\infty}}_k(q,n,d):=    \max \{ \sum_{v \in \Z_q^n} x(\{v\})\,\, &|\,\,x:\mathcal{C}_k \to \R, \,\, x(\emptyset )=1,\,\, x(S)=0 \text{ if~$d_{\text{min}}^{L_{\infty}}(S)<d$}, \notag \\[-1.1em]
& \,\, M_{k,D}(x) \text{ is positive semidefinite for each~$D$ in~$\mathcal{C}_k$}\}. 
\end{align}
So~$B_k^{L_{\infty}}(q,n,d)$ is obtained from the bound~$B_k^{L}(q,n,d)$ in~$(\ref{Bleend})$ by replacing in the definition $d_{\text{min}}^{L}(S)$ by~$d_{\text{min}}^{L_{\infty}}(S)$. It is not hard to see that~$\alpha(C_{d,q}^n) \leq B_k^{L_{\infty}}(q,n,d)$, by a proof analogous to that of Proposition~\ref{trivleeprop}. For comparison,~$\vartheta(C_{d,q}^n)$ is equal to the bound obtained from~$B_2^{L_{\infty}}(q,n,d)$  by removing the constraints that~$M_{2,D}(x)$ is positive semidefinite for subsets~$D \in \mathcal{C}_2$ with~$D \neq \emptyset$. Moreover,~$B_2^{L_{\infty}}(q,n,d)$ is equal to the Delsarte bound, which is equal to the bound~$\vartheta'(C_{d,q}^n)$, with~$\vartheta'$ as in~$\cite{thetaprime}$.

\begin{table}[ht]
\centering
\begin{tabular}[t]{l|ccccc}
 & 1 & 2 & 3 & 4 & 5 \\  \hline
$B_2^{L_{\infty}}(5,n,2)$ &  $2.236 $ & $ 5.000 $ & 11.180  & $25.000$ & $55.902$ \\
$B_3^{L_{\infty}}(5,n,2)$ & 2.000 & $5.000$  & $10.915$ & $25.000$ & $55.902$ \\ \hline
$B_2^{L_{\infty}}(7,n,2)$ & $3.318$ & $ 11.007 $ & 36.517  & $121.152$ & $401.943$ \\
$B_3^{L_{\infty}}(7,n,2)$ & 3.000 & $10.260$  & $35.128$ & $119.537$ & $401.908$\\ \hline
$B_2^{L_{\infty}}(7,n,3)$ & $2.110$ & $ 4.452 $ & $9.393 $  & $19.818$ &  $ 41.814$    \\
$B_3^{L_{\infty}}(7,n,3)$ & 2.000 & $4.139$  & $8.957$ & $19.494$ & $41.782 $\\ \hline
$\#$Vars in $B_3^{L_{\infty}}(5,n,2)$ & 2 & 9  & 48 & 214 & 799 \\
$\#$Vars in $B_3^{L_{\infty}}(7,n,2)$ & 3 & 43  & 423 & 3161 & 19023 \\
$\#$Vars in $B_3^{L_{\infty}}(7,n,3)$ & 2 & 12  & 137 & 1316 & 9745 \\
\end{tabular}
\caption{Bounds on~$\alpha(C_5^n)$,~$\alpha(C_7^n)$ and~$\alpha(C_{3,7}^n)$, rounded to three decimal places. It holds that $B_2^{L_{\infty}}(5,n,2) = \sqrt{5}^n$.\label{b3maxtable} }
\end{table}

To compute~$B_3^{L_{\infty}}(q,n,d)$, the reductions from Section~\ref{reduct} can be used. The new bound $B_3^{L_{\infty}}(q,n,d)$ does not seem to improve significantly over the bound obtained from Lov\'asz's $\vartheta$-function, except for very small~$n$. See Table~\ref{b3maxtable} for some results for~$q \in \{5,7\}$ and~$1 \leq n \leq 5$. For these cases,~$B_3^{L_{\infty}}(q,n,d)$ does not give new upper bounds on~$\alpha(C_{d,q}^n)$, as the values~$\alpha(C_5^3)=10$, $\alpha(C_7^2)=10$,~$\alpha(C_7^3)=33$ (cf.~\cite{baumert}),~$\alpha(C_{3,7}^3)=8$ (cf.~\cite{circular}) are already known and~$\alpha(C_7^4)\leq \floor{(7/2)\alpha(C_7^3)}=115$. The number of variables ``$\#$Vars'' in~$B_3^{L_{\infty}}(q,n,d)$, which  is the number of~$D_q^n \rtimes S_n$-orbits of nonempty codes of size~$\leq 3$ and minimum Lee\textsubscript{$\infty$}-distance at least~$d$,  is also given in Table~\ref{b3maxtable} for the considered cases.

\appendices

\section{Appendix: Formulas with integers\label{intformulaslee}}  \small

Note that for~$q \leq 4$ or~$q=6$ all coefficients the formulas in Section~\ref{d0comp} for the block diagonalization of~$M_{3,\emptyset}(z)$ are rational (hence all constraints can be made integer). For other~$q $ the formulas contain irrational numbers. To obtain a semidefinite program which only contains integers we used in the implementation for~$q=5$ and~$q=7$ not the representative set from~$(\ref{reprsetdq2})$ for the action of~$D_q^n \rtimes S_n$ on~$\C^{\Z_q^n}$  but the representative set from~$(\ref{reprsetdq1})$ for the action of~$S_2^n \rtimes S_n$  on~$\C^{\Z_q^n}$  to reduce the matrix~$\sum_{\omega \in \Omega_2} z(\omega) N_{\omega}'$. Then~$(\ref{PhiR})$ gives (where we write~$\Psi$ for the map in~$(\ref{PhiR})$  to distinguish it from the map~$\Phi$ from Section~\ref{d0comp})
\begin{align} \label{blocks1ttdpsi}
    \Psi \left(\sum_{\omega \in \Omega_2}z(\omega) N_{\omega}' \right ) = \bigoplus_{\bm{n} \in \bm{N}} \bigoplus_{\bm{\lambda} \vdash \bm{n}} \sum_{\omega \in \Omega_2}z(\omega) U_{\bm{\lambda}}^T N_{\omega}' U_{\bm{\lambda}},
\end{align}
where~$U_{\bm{\lambda}}$ denotes the matrix in~$(\ref{reprsetdq1})$ that corresponds with~$\bm{\lambda} \in W_{\bm{\lambda}}$. By~$\eqref{kmubelangrijk}$, we have $\sum_{\mu \in Q} v_{\bm{\tau}}^T K_{\mu} v_{\bm{\sigma}} \mu =p_{\bm{\tau,\sigma}}$. From this, one obtains that $\sum_{\mu \in Q'} v_{\bm{\tau}}^T K_{\mu}' v_{\bm{\sigma}} \mu =p_{\bm{\tau,\sigma}}'$, where~$p_{\bm{\tau,\sigma}}'$ is the polynomial obtained from~$p_{\bm{\tau,\sigma}}$ by replacing each variable~$d_{0ij}^* \in Z^*$ (with~$0ij\in \Pi$) by the variable~$f_{\pi(iij)}^* \in \widetilde{Z}^*$. Hence the following replacements must be done, using the formulas from~\eqref{formuleslee1}:
\begin{alignat}{2} 
B_1(1) \otimes B_1(1) &= 1 d_{000}^*  &&\mapsto 1 f_{000}^*, \notag   \\ 
B_1(1) \otimes B_1(j+1) &= 2 d_{00j}^* &&\mapsto 2 f_{00j}^*,    \text{ for~$j=1,\ldots,\floor{q/2}$} \notag \\
B_1(j+1) \otimes B_1(1) &=  2 d_{0j0}^* &&\mapsto 2 f_{00j}^*,    \text{ for~$j=1,\ldots,\floor{q/2}$}  \notag\\
B_1(j+1) \otimes B_1(h+1) &=2  d_{0jh}^* + 2 d_{0j(q-h)}^* &&\mapsto 2 f_{00t_1}^* + 2 f_{00t_2}^*,    \text{ for~$j,h \in \{1,\ldots,\floor{q/2}\}$}, \notag  \\ 
B_2(j) \otimes B_2(h) &= 2 d_{0jh}^* - 2 d_{0j(q-h)}^* &&\mapsto 2 f_{00t_1}^* - 2 f_{00t_2}^*,    \text{ for~$j,h \in \{1,\ldots,\floor{(q-1)/2}\}$},   \label{formuleslee2}
\end{alignat}  
where in the above formulas we set~$t_1 := j-h $ if~$j \geq h$ and~$t_1:= h-j$ else, % so that~$t_1 \in \{0,\ldots,\floor{q/2}-1\}$,
  and we set~$t_2:= j+h$ if~$j+h\leq \floor{q/2}$ and~$t_2:=q-(j+h)$ else.
  
So one computes the entry $\sum_{\omega \in \Omega_2} z(\omega) \bm{v_{\bm{\tau}}}^T   N_{\omega}' \bm{v_{\bm{\sigma}}}$ in the block $ \sum_{\omega \in \Omega_2}z(\omega) U_{\bm{\lambda}}^T N_{\omega}' U_{\bm{\lambda}}$ by first expressing~$p_{\bm{\tau}, \bm{\sigma}}'$ as a linear combination of~$\mu \in Q'$ and subsequently replacing each $\mu \in Q'$ in~$p_{\bm{\tau}, \bm{\sigma}}'$ with the variable~$z(r'(\mu))$.

To add the empty code, one may add a new unit base vector~$\epsilon_{\emptyset}$ to the matrix in the representative set~$\eqref{reprsetdq1}$ indexed by~$\bm{n}=((n),())$ and calculate the new entries~$ \epsilon_{\emptyset}^T M_{2,\emptyset}(z)   v_{\bm{\sigma}}$, for each~$\sigma \in W_{\bm{\lambda}}$. However, this is not necessary. As
$\Psi \left(\sum_{\omega \in \Omega_2}z(\omega) N_{\omega}' \right )$ is positive semidefinite if and only if~$\Phi \left(\sum_{\omega \in \Omega_2}z(\omega) N_{\omega}' \right )$  is positive semidefinite,
 and~$M_{2,\emptyset}(z)$ is positive semidefinite if and only if both~$\Phi \left(\sum_{\omega \in \Omega_2}z(\omega) N_{\omega}' \right )$  and~$T$ from~$\eqref{emptycodeblock}$ are positive semidefinite, we find that 
 $$
M_{2,\emptyset}(z) \text{ is positive semidefinite } \Longleftrightarrow  \Psi \left(\sum_{\omega \in \Omega_2}z(\omega) N_{\omega}' \right ) \text{ and } T \text{ are positive semidefinite}.
 $$  
So the~$2\times 2 $ matrix~$T$ together with the matrix blocks in~$\eqref{blocks1ttdpsi}$ form a block diagonalization of~$M_{2,\emptyset}(z)$.\iffalse\footnote{An advantage of this approach  is that it makes it possible to  exploit the fact that the second row and column of~$T$ are divisible by~$q^n$. One can divide the second row and column of~$T$ by~$q^n$ (so the bottom right entry by~$q^{2n}$) and also the objective function by~$q^n$, to obtain a more stable semidefinite program.}\fi

\section{Appendix: An overview of the program}\small
In this section we give a high-level overview of the program, to help the reader with implementing the method. 
See Figure~$\ref{pseudocode}$ for an outline of the method.

A few remarks regarding the implementation:
\begin{enumerate}[(i)]
\item We write~$\omega_0$ for the unique~$D_q^n\rtimes S_n$-orbit corresponding to a code of size~$1$.
\item To speed up the replacement of monomials in~$d_P^*$ or~$f_P^*$ by variables~$z(\omega)$, it is useful to add a preprocessing step to determine in advance for each degree~$n$ monomial~$\mu = d_{P_1}^*\ldots d_{P_n}^*$ with all~$P_i \in \Pi$ and~$\mu' = f_{P_1}^*\ldots f_{P_n}^*$ with all~$P_i \in \Pi'$ which orbit~$r(\mu) \in \Omega_3$ or~$r'(\mu') \in \Omega_2$ corresponds with it. If the orbit corresponds to a code of minimum Lee (or Lee\textsubscript{$\infty$})  distance~$<d$ to zero, we must set the corresponding variable to zero and can delete it from the program. 

\begin{figure}[H]
  \fbox{
    \begin{minipage}{15.25cm}\small
     \begin{tabular}{l}
    \textbf{Input: } Natural numbers~$q,n$ \\
	\textbf{Output: }Semidefinite program to compute~$B_3^L(q,n,1)$ $(=B_3^{L_{\infty}}(q,n,1))$.  \\
	$\,$\\
    \printv \emph{Maximize} $q^n z(\omega_0)$  \\
        \printv \emph{Subject to:}  \\    
\text{\color{blue}//Start with~$|D|=1$.}\\
\foreachv $\bm{n}=(n_1,n_2) \in \bm{N}$  \\
	\hphantom{1cm} \foreachv $\bm{\lambda}=(\lambda_1,\lambda_2) \vdash \bm{n}$ with height$(\lambda_1) \leq \floor{q/2}+1 $,   height$(\lambda_2) \leq \floor{(q-1)/2} $\\
			\hphantom{1cm} \hphantom{1cm} start a new block~$M_{\bm{\lambda}}$\\
			\hphantom{1cm} \hphantom{1cm} \foreachv $\bm{\tau} \in W_{\bm{\lambda}}$ from~$(\ref{wlambda})$\\
					\hphantom{1cm} \hphantom{1cm}\hphantom{1cm} \foreachv $\bm{\sigma} \in W_{\bm{\lambda}}$ from~$(\ref{wlambda})$\\
					\hphantom{1cm} \hphantom{1cm}\hphantom{1cm}\hphantom{1cm} compute~$p_{\bm{\tau},\bm{\sigma}}$ from~$(\ref{polylee})$ as linear combination in~$B_i(j)\otimes B_i(h)$\\			
					\hphantom{1cm} \hphantom{1cm}\hphantom{1cm}\hphantom{1cm} replace each~$B_i(j)\otimes B_i(h)$ by the linear expression in~$d_P^*$ from~\eqref{formuleslee1} \\	
							\hphantom{1cm} \hphantom{1cm}\hphantom{1cm}\hphantom{1cm}		replace each degree~$n$ monomial~$\mu$ in~$d_{P}^*$ by a variable~$z(r(\mu))$\\
	\hphantom{1cm} \hphantom{1cm}\hphantom{1cm}\hphantom{1cm}	 $(M_{\bm{\lambda}})_{\bm{\tau},\bm{\sigma}}:= $ the resulting linear polynomial in variables~$z(\omega)$ \\				
					\hphantom{1cm} \hphantom{1cm}\hphantom{1cm} \ndv\\
			\hphantom{1cm} \hphantom{1cm} \ndv\\
		  	\hphantom{1cm} \hphantom{1cm}   \printv $M_{\bm{\lambda}}$\emph{ positive semidefinite}.   \\
			\hphantom{1cm} \ndv\\  	
 \ndv\\
 \text{\color{blue}//Now~$D=\emptyset$.}\\
\foreachv $\bm{n} \in \bm{M}$  \\
	 \hphantom{1cm} start a new $(1 \times 1)$-block~$M_{\bm{n}}$  \\
      \hphantom{1cm} compute~$p_{\bm{n}}$ from~$(\ref{polylee2})$ as linear combination in~$C_i\otimes C_i$\\			
					\hphantom{1cm}  replace each~$C_i\otimes C_i$ by the linear expression in~$f_P^*$ from~\eqref{formuleslee2compact} \\	
							\hphantom{1cm}		replace each degree~$n$ monomial~$\mu$ in~$f_{P}^*$ by a variable~$z(r'(\mu))$\\
	                           \hphantom{1cm}	 $(M_{\bm{n}}):= $ (the resulting linear polynomial in variables~$z(\omega)$) \\	
	                           	                           	\hphantom{1cm}   \ifv $\bm{n}=(n,0,\ldots,0)$ \,\,\,\,\,   \text{\color{blue}//add a row and a column corresponding to $\emptyset$.}\\		
		                          \hphantom{1cm} 	\hphantom{1cm} add a row and column to $M_{\bm{n}}$ indexed by~$\emptyset$   \\		   
		                         \hphantom{1cm}  	\hphantom{1cm} put $(M_{\bm{n}})_{\emptyset,\emptyset}:=1$ and $(M_{\bm{n}})_{\bm{n},\emptyset} =(M_{\bm{n}})_{\emptyset,\bm{n}}:=q^n z_{\omega_0}$  \\	                        	                           	
	                           	\hphantom{1cm} \ndv   \\			                           	                           	
	                           	\hphantom{1cm}   \printv $M_{\bm{n}}$\emph{ positive semidefinite}.   \\			
					\ndv\\
	\text{\color{blue}//Now nonnegativity of all variables.}\\
		\foreachv $\omega \in \Omega_3$  \\
	    \hphantom{1cm} \printv $ z(\omega) \geq 0$\\
	    \ndv 
     \end{tabular}
    \end{minipage}    } \caption{\label{pseudocode}\small{Algorithm to generate a semidefinite program for computing~$B_3^L(q,n,1)$. To compute~$B_3^L(q,n,d)$ or~$B_3^{L_{\infty}}(q,n,d)$, one must set all variables~$z(\omega)$  with~$\omega \in \Omega_3$ an orbit corresponding to a code of minimum Lee (respectively, Lee\textsubscript{$\infty$}) distance~$<d$ to zero. 
If rows and columns in matrix blocks~$M_{\bm{\lambda}}$ consist only of zeros after the replacement, it is useful to remove these rows and columns.}}
\end{figure} \small

\item In case~$D= \emptyset$, the matrix blocks contain irrational numbers for~$q \notin \{2,3,4,6\}$. In Appendix~\ref{intformulaslee} it is explained how to obtain a semidefinite program which only contains integers. This is not displayed in the above pseudocode, but the adaptations are straightforward.

 In order to obtain the matrix blocks from~$\eqref{blocks1ttdpsi}$ for~$D=\emptyset$ one can simply repeat the steps in the above pseudocode for~$|D|=1$, but with the following adaptation: replace each~$B_i(j) \otimes B_i(h)$ by the linear expression in~$f_P^*$ from~$(\ref{formuleslee2})$, and subsequently replace each monomial~$\mu$ of degree~$n$ in~$f_P^*$ by a variable~$z(r'(\mu))$.

\item The programs we used to generate input for the SDP-solver can be found at the following location (also accessable via the author's website):
$$
\text{\url{https://drive.google.com/open?id=1-XRbfc4TYhoySC33GRWfvNEMOZEltg6X}}. 
$$

\end{enumerate}

 \section*{Acknowledgements}
 The author wants to thank Lex Schrijver for very useful discussions. Also, the author wants to thank Sander Gribling, Bart Litjens, Bart Sevenster and the anonymous referees for useful comments.  The research leading to these
results has received funding from the European Research Council under the European Union’s Seventh Framework Programme (FP7/2007-2013) / ERC grant agreement \textnumero 339109.
\footnotesize
\selectlanguage{english} 

\end{document}